\documentclass[11pt]{article}
\usepackage{amssymb,amsmath}
\usepackage{amsmath}
\usepackage{amsfonts}
\newtheorem{precor}{{\bf Corollary}}

\newenvironment{cor}{\begin{precor}{\hspace{-0.5
               em}{\bf.\ }}}{\end{precor}}
\newtheorem{precon}{{\bf Conjecture}}

\newenvironment{con}{\begin{precon}{\hspace{-0.5
               em}{\bf.\ }}}{\end{precon}}
\newtheorem{predefin}{{\bf Definition}}

\newtheorem{preexm}{{\bf Example}}

\newtheorem{preappl}{{\bf Application}}

\newtheorem{prelem}{{\bf Lemma}}

\newenvironment{lem}{\begin{prelem}{\hspace{-0.5
               em}{\bf.\ }}}{\end{prelem}}
\newtheorem{preproof}{{\bf Proof.\ }}

\newenvironment{proof}[1]{\begin{preproof}{\rm
               #1}\hfill{$\blacksquare$}}{\end{preproof}}
\newtheorem{presproof}{{\bf Sketch of Proof.\ }}

\newtheorem{prethm}{{\bf Theorem}}

\newenvironment{thm}{\begin{prethm}{\hspace{-0.5
               em}{\bf.\ }}}{\end{prethm}}
\newtheorem{prealphthm}{{\bf Theorem}}

\newtheorem{prealphlemma}{{\bf Lemma}}

\newtheorem{prepro}{{\bf Proposition}}

\newtheorem{preprb}{{\bf Problem}}

\newtheorem{prequ}{{\bf Question}}

\def\conct[#1,#2]{\mbox {${#1} \leftrightarrow {#2}$}}
\def\dconct[#1,#2]{\mbox {${#1} \rightarrow {#2}$}}
\def\deg[#1,#2]{\mbox {$d_{_{#1}}(#2)$}}
\def\mindeg[#1]{\mbox {$\delta_{_{#1}}$}}
\def\maxdeg[#1]{\mbox {$\Delta_{_{#1}}$}}
\def\outdeg[#1,#2]{\mbox {$d_{_{#1}}^{^+}(#2)$}}
\def\minoutdeg[#1]{\mbox {$\delta_{_{#1}}^{^+}$}}
\def\maxoutdeg[#1]{\mbox {$\Delta_{_{#1}}^{^+}$}}
\def\indeg[#1,#2]{\mbox {$d_{_{#1}}^{^-}(#2)$}}
\def\minindeg[#1]{\mbox {$\delta_{_{#1}}^{^-}$}}
\def\maxindeg[#1]{\mbox {$\Delta_{_{#1}}^{^-}$}}
\def\isdef{\mbox {$\ \stackrel{\rm def}{=} \ $}}
\def\dre[#1,#2,#3]{\mbox {${\cal E}_{_{#3}}(#1,#2)$}}
\def\pdre[#1,#2,#3]{\mbox {${\cal P}_{_{#3}}(#1,#2)$}}
\def\var[#1,#2]{\mbox {${\rm Var}_{_{#1}}(#2)$}}
\def\ls[#1]{\mbox {$\xi^{^{#1}}$}}
\def\hom[#1,#2]{\mbox {${\rm Hom}({#1},{#2})$}}
\def\onvhom[#1,#2]{\mbox {${\rm Hom^{v}}(#1,#2)$}}
\def\onehom[#1,#2]{\mbox {${\rm Hom^{e}}(#1,#2)$}}
\def\core[#1]{\mbox {$#1^{^{\bullet}}$}}
\def\cay[#1,#2]{\mbox {${\rm Cay}({#1},{#2})$}}
\def\cays[#1,#2]{\mbox {${\rm Cay_{s}}({#1},{#2})$}}
\def\dirc[#1]{\mbox {$\stackrel{\rightarrow}{C}_{_{#1}}$}}
\def\cycl[#1]{\mbox {${\bf Z}_{_{#1}}$}}

\def\sdg[#1]{\mbox {$\stackrel{\leftrightarrow}{#1}$}}
\begin{document}
\begin{center}
{\Large \bf On Dynamic Coloring of Graphs}\\
\vspace*{0.5cm}
{\bf Meysam Alishahi}\\
{\it Department of Mathematical Sciences}\\
{\it Shahid Beheshti University, G.C.}\\
{\it P.O. Box {\rm 1983963113}, Tehran, Iran}\\
{\tt m\_alishahi@sbu.ac.ir}\\

\end{center}
\begin{abstract}
\noindent A dynamic coloring of a graph $G$ is a proper coloring
 such that for every vertex $v\in V(G)$ of degree at least 2,
the neighbors of $v$ receive at least 2 colors. In this paper we
present some upper bounds for the dynamic chromatic number of
graphs. In this regard, we shall show that there is a constant $c$ such
that for every $k$-regular graph $G$, $\chi_d(G)\leq \chi(G)+
c\ln k +1$. Also, we introduce an upper bound for the dynamic
list chromatic number of regular graphs.
\begin{itemize}
\item[]{{\footnotesize {\bf Key words:}\  Dynamic chromatic number, Dynamic list chromatic number.}}
\item[]{ {\footnotesize {\bf Subject classification: 05C} .}}
\end{itemize}
\end{abstract}
\section{Introduction}
Let $H$ be a hypergraph. The vertex set and the
hyperedge set of $H$  are mentioned as $V(H)$ and
$E(H)$, respectively. The maximum degree and the minimum
degree of $H$ are denoted by $\Delta(H)$ and
$\delta(H)$, respectively.
For an integer $l\geq 1$  we denote by $[l]$ the set $\{1, 2,
\ldots, l\}$. A proper $l$-coloring of a hypergraph
$H$ is a function $c: V(H)\longrightarrow [l]$  in which there is no monochromatic
hyperedge in $H$. We say a hypergraph $H$ is
$t$-colorable if there is a proper $t$-coloring of it.  For a
hypergraph $H$, the smallest integer $l$ that
$H$ is $l$-colorable is called the chromatic number of
 $H$ and denoted by $\chi(H)$. Note that a graph $G$ is a
hypergraph such that  the cardinality of each $e\in E(G)$ is 2. We say
a graph $G$ is $t$-critical if $\chi(G)=t$ and any proper induced
subgraph of $G$ has the chromatic number strictly less than $t$.
For a vertex $v\in V(G)$, $N(v)$ is the set of all adjacent
vertices to $v$ and ${\rm deg}_T(v)$ is the number of neighbors of
$v$ that are lied in $T\subseteq V(G)$, i.e., the cardinality of $N(v)\cap T$.

We denote by ${[m] \choose n}$ the collection of all $n$-subsets
of $[m]$. The {\it Kneser graph} $KG(m,n)$ is the graph with
vertex set ${[m] \choose n}$, in which $A$ is connected to $B$ if
and only if $A \cap B = \emptyset$. It was conjectured by Kneser
\cite{MR00685361} in 1955, and proved by Lov\'{a}sz
\cite{MR514625} in 1978, that $\chi(KG(m,n))=m-2n+2$.

A proper vertex $l$-coloring of a graph $G$ is called a dynamic
$l$-coloring \cite{Montgomery} if for every vertex $u$ of degree
at least $2$, there are at least two different colors appeared in
the neighborhood of $v$. The smallest integer $l$ that there is
a dynamic $l$-coloring of $G$ is called {\it the dynamic
chromatic number of $G$} and denoted by $\chi_d(G)$. Obviously,
$\chi(G)\leq\chi_d(G)$. Some properties of dynamic coloring was
studied in \cite{akbari2,MR2251583,MR1991048,Montgomery,123}.
It was proved in \cite{MR1991048} that for a connected graph $G$ if  $\Delta\leq3$,
then $\chi_d(G)\leq4$ unless $G=C_5$, in which case
$\chi_d(C_5)=5$ and if $\Delta\geq4$, then $\chi_d(G)\leq\Delta+1$.
It was shown in \cite{Montgomery} that the difference between chromatic number and
dynamic chromatic number can be arbitrarily large. However, it was
conjectured that for regular graphs the difference is at most 2.
\begin{con}{\rm\cite{Montgomery}}\label{conj}
For any regular graph $G$, $\chi_d(G)-\chi(G)\leq 2$
\end{con}
Also, it was proved in \cite{Montgomery} that if $G$ is a
bipartite $k$-regular graph, $k\geq 3$ and $n<2^k$ then
$\chi_d(G)\leq 4$.

For a graph $G$, $L$ is called an $l$-list assignment for $G$
if for each vertex $v\in V(G)$, $L(v)$ is an $l$-set
of available colors  at $v$. An $L$-coloring, is a proper coloring
$c$ such that $c(v)\in L(v)$, for each $v \in V(G)$. The graph G
is $l$-list colorable if for every $l$-list assignment $L$ to the
vertices of $G$, $G$ has  a proper $L$-list coloring. The list
chromatic number, $\chi_l(G)$ is the minimum number $l$ such that
G is $l$-list colorable. The {\it list dynamic chromatic number}
of a graph $G$, ${\rm ch}_d(G)$, is the minimum positive integer
$l$ such that for every $l$-list assignment, there is a dynamic
coloring of $G$ such that every vertex of $G$ is colored with a
color from its list. Clearly, the dynamic list chromatic number of graphs is a
common generalization of both dynamic chromatic and list chromatic
number of graphs and also $\chi_d(G)\leq {\rm ch}_d(G)$. It was proved in
\cite{akbari2} that for a connected graph $G$ if  $\Delta(G)\geq3$
then ${\rm ch}_d(G)\leq \Delta(G)+1$, and if $\Delta(G)\leq 3$
then ${\rm ch}_d(G)\leq 4$ except $G=C_5$, in which case ${ch_2}(C_5)=5$.

In a graph $G$, a set $T\subseteq V(G)$ is called a {\it total
dominating set} in $G$ if for every vertex $v \in V(G)$, there is
at least one vertex $u \in T$ adjacent to $v$. The set $T$ is
called {\it double total dominating set} if both of $T$ and its
complement $V(G)\setminus T$ are total dominating set.

\section{Results}
Suppose every vertex of a graph $G$ appears in some
triangles. It is clearly that $\chi_d(G)=\chi(G)$. We shall use
this simple result to prove the next theorem.
\begin{thm}
Let $0<p<1$ be a constant. Almost all graphs in ${\cal G}(n,p)$
have the same chromatic and dynamic chromatic number.
\end{thm}
\begin{proof}{
We show that for almost all graphs $G\in{\cal G}(n,p)$, every
vertex $v\in V(G)$ appears in some triangles. Consider complete graph $K_n$ with
$V(K_n)=V({\cal G}(n,p))=[n]$. Let $v$ be an arbitrary
vertex in $V(G)$. Consider $\lfloor\frac{n-1}{2} \rfloor$ edge disjoint triangles in $K_n$ such
that they all have the vertex $v$. Note that no vertex, except $v$, is in
more than one triangle. Define ${\cal A}_v$ to be the event that
none of these triangles is happened in $G$. Clearly, ${\rm pr}({\cal
A}_v)\leq (1-p^3)^{\frac{n-2}{2}}$. As
$n(1-p^3)^{\frac{n-2}{2}}\longrightarrow 0$æ the proof is
completed. }\end{proof}

In the next theorem, we present an upper bound for the dynamic
chromatic number of $k$-regular graph $G$ in terms of
$k$ and $\chi(G)$.
\begin{thm}\label{main}
There exists a constant $c$ such that for any $k$-regular graph
$G$, $\chi_d(G)\leq \chi(G)+ c\ln k +1$.
\end{thm}
\begin{proof}{
Define $p=c'{\ln k\over k}$ where $c'$ will be specified later.
Consider a random set  $T$ such that each vertex $u$ lies in $T$
with the probability $p$. Assume that for each vertex $u$, the
random variable $X_u$ is the number of neighbors of $u$ that are
in $T$. Clearly, $X_u$ is a binomial random variable and according to the
Chernoff inequality we have
$${\rm Pr}(|X_u-E(X_u)|> \lambda)\leq2e^{-{\lambda^2\over 3E(X_u)}}, \hspace*{2cm}  0<\lambda<E(X_u)$$
Obviously, $E(X_u)=c'\ln k$. For each vertex $u$ define ${\cal
A}_u$ to be the event that $|{\rm deg}_T(u)-c'\ln k|> \lambda$. By
Chernoff inequality we have ${\rm Pr({\cal A}_u)}\leq
2e^{-{\lambda^2\over 3E(X_u)}}.$ For each $u$, ${\cal A}_u$ is
mutually independent of all ${\cal A}_v$ events but at most $k^2$
number of them. If we set $c'>6$ then there exist a threshold
$n(c')$ such that when $k\geq n(c')$ we have
$$3c'\ln k(1+\ln 2+\ln(k^2+1))<{c'}^2(\ln k)^2.$$
Let $k\geq n(c')$ and consider a $\lambda$ such that
\begin{equation}\label{equ}
3c'\ln k(1+\ln 2+\ln(k^2+1))\leq\lambda^2<{c'}^2(\ln k)^2.
\end{equation}
In view of previous inequality, we have
$$2e(k^2+1)e^{-{\lambda^2\over 3E(X_u)}}\leq 1.$$
Therefore, By applying Lovasz Local Lemma there exists a set
$T\subseteq V(G)$ such that for every $v\in V(G)$, ${\cal A}_v$
does not happen. Equivalently, there is a set $T\subseteq V(G)$
such that for every vertex $v\in V(G)$, $|{\rm deg}_T(v)-c'\ln
k|\leq \lambda$. Consequently, since $0<\lambda< c'\ln k$, for
each vertex $v\in V(G)$ we have $0<{\rm deg}_T(v)<2c'\ln k$. It
implies that for every vertex $u\in V(G)$, ${\rm deg}_T(u)>0$ and
${\rm deg}_{V(G)\setminus T}(u)>0$. Note that $\Delta (G[T])\leq
2c'\ln k$ and therefore $\chi(G[T])=l\leq \lfloor c\ln
k\rfloor+1$. Color the vertices in $T$ with colors come from
$\{1,2,\ldots,l\}$ and also color the vertices in $V(G)\setminus
T$ with colors come from $\{l+1,l+2,\ldots,l+\chi(G)\}$. One can
see that this coloring is a dynamic coloring of $G$ and uses at
most $\chi(G)+\lfloor c\ln k\rfloor+1$ colors. Since the number
of $k< n(c')$ is finite, there is a constant $c$ such that for
every $k$-regular graph $G$, $\chi_d(G)\leq \chi(G)+ c\ln k +1$
and the proof is completed. }\end{proof}

It was conjectured in \cite{akbari2} that for any graph $G$,
${\rm ch}_d(G)=\max \{\chi_l(G), \chi_d(G)\}$. This conjecture was disproved in
\cite{1234}. It was shown in \cite{1234} that for any integer $k \geq 5$,
there is a bipartite graph $G_k$ such that ${\rm ch}(G_k) = \chi_d(G_k)=3$ and ${\rm ch}_d(G)\geq k$.

The next theorem provides an upper bound for list dynamic chromatic
number of regular graphs in terms of their list chromatic number.
\begin{thm}
Let  $\epsilon$ be a positive constant. If $G$ is a $k$-regular
graph then for large enough $k$, $\chi_d(G)\leq {\rm ch}_d(G)\leq \lceil
(1+\epsilon)\chi_l(G)\rceil$
\end{thm}
\begin{proof}{
For convenience let $\chi_l(G)=l$. Without loss of generality, we
can assume that $G$ is not a complete graph. Let $m\geq l$ be a
positive integer and consider an $m$-list assignment $L$ for $G$.
For each vertex $v\in V(G)$, choose a random
$l$-set $L'(v)\subset L(v)$ uniformly and independently. For $v\in
V(G)$, suppose that ${\cal B}_v$ denotes the event that
$\cap_{u\in N(v)}L(u)\neq \varnothing$. It is readily seen that
${\rm Pr}({\cal B}_v)\leq m({l\over m})^k$. Also, ${\cal B}_v$ is
mutually independent of all the other events  ${\cal B}_u$ but
those for which $N(v)\cap N(u)\neq\varnothing$. Therefore, at most
$k(k-1)$ events are not mutually independent of ${\cal B}_v$. In
view of Lovasz Local Lemma, if $m$ is sufficiently large such that
\begin{equation}\label{eq}
ek^2m({l\over m})^k\leq 1
\end{equation}
 then there is a list assignment $L'$
such that for each vertex $v\in V(G)$, $L'(v)\subseteq L(v)$,
$|L'(v)|=l$ and $\cap_{u\in N(v)}L'(u)=\varnothing$. Obviously,
$ek^2m({l\over m})^k\leq 1$ if and only if
$l(elk^2)^{\frac1{k-1}}\leq m$. Note that $l\leq k+1$ and
therefore $(elk^2)^{\frac1{k-1}}\longrightarrow 1$. Hence, there
is threshold $M(\epsilon)$ such that if $k\geq M(\epsilon)$ then
$(elk^2)^{\frac1{k-1}}< 1+\epsilon$. Hence, if $k\geq
M(\epsilon)$ then $m=\lceil (1+\epsilon)l\rceil$ satisfies
Equation \ref{eq}. Since, for every $v\in V(G)$,
$|L'(v)|=l=\chi_l(G)$, $G$ has an $L$-coloring $c$ such
that $c(u)\in L'(u)\subseteq L(u)$ for each $u\in V(G)$. Note that for every
vertex $v\in V(G)$, $\cap_{u\in N(v)}L'(u)=\varnothing$ and it
obviously implies that $c$ is a dynamic coloring of $G$.
}\end{proof}
Let $H$ be a hypergraph. 2-colorability of hypergraphs has been
studied in the literature and has lots of applications in some other
concepts of combinatorics. Hereafter, we want to make a
connection between 2-colorability of  hypergraphs and
the dynamic chromatic number of graphs.
\begin{lem}\label{hyper}
Let H be a hypergraph. Assume that there are two positive integers
$2\leq r_1\leq r_2$ such that for any $f\in E(H)$, $
r_1\leq|f|\leq r_2$ and also
$er_2\Delta(H)({1\over2})^{r_1-1}\leq1$. Then $H$ is
$2$-colorable.
\end{lem}
\begin{proof}{
Consider a random 2-coloring of $V(H)$ in which each vertex $v\in
V(H)$ is colored red or blue with the same probabilities. For any
$f\in E(H)$, define ${\cal A}_f$ to be the event that {\it all
the vertices in $f$ have the same color}. Obviously, ${\rm
Pr}({\cal A}_f)\leq({1\over2})^{r_1-1}$ and two events ${\cal
A}_f$ and ${\cal A}_f'$ are mutually independent when $f\cap
f'=\varnothing$. So, there are at most $r_2(\Delta (H)-1)+1$
events that are not mutually independent of ${\cal A}_f$. Lovasz
Local Lemma implies that if
$er_2\Delta(H)({1\over2})^{r_1-1}\leq1$, then $H$ is
$2$-colorable. }\end{proof}

It was shown in \cite{Montgomery} that for any $l\geq 2$ there
is a family of $l$-colorable graphs such that the difference between chromatic and dynamic chromatic
number, in this family, is unbounded. The next theorem states that if $\Delta(G)$ and $\delta(G)$
are not so far from each other then its chromatic
and  dynamic chromatic number are not far
from each other.
\begin{lem}\label{delta}
Let $G$ be a graph such that  $e\Delta^2\leq2^{\delta(G)-1}$. Then $\chi_d(G)\leq 2\chi(G).$
\end{lem}
\begin{proof}{
Define a hyprgraph $H$ whose vertex set is the same as the vertex
set of $G$ and its  hyperedge set is defined as follows. $$E(H)\isdef
\{N(v)|v\in V(G)\}.$$
Clearly, for any $f\in E(H)$, $\delta(G)\leq |f|\leq \Delta(G)$
and $\Delta(G)=\Delta(H)$. By considering Lemma \ref{hyper}, $H$
is 2-colorable. Let $f$ be a 2-coloring of $H$ and $c$ be a
$\chi(G)$-coloring of $G$. It is readily to seen $h=(f,c)$ is a
$2\chi(G)$ dynamic coloring of $G$. }\end{proof}

Note that when $G$ is a $k$-regular graph and $k\geq 9$, Lemma
\ref{delta} implies that $\chi_d(G)\leq 2\chi(G)$. It was shown
by Thomassen \cite{thomassen} that for any $k$-uniform and
$k$-regular hypergraph $H$, if $k\geq 4$ then $H$ is 2-colorable.

Assume that $H$ is a $k$-uniform hypergraph and $\Delta(H)\leq k$.
One can  construct a hypergraph $H'\supseteq H$ such that $H'$ is
$k$-uniform and $k$-regular as follows. If $\delta(H)=k$ then
$H'=H$. Assume that $H^i$ is constructed and $k-\delta(H^i)=t>0$.
Consider $H^{*}$ to be a hypegraph which is a union of $k$ disjoint
copies of $H^{i}$. For each vertex $v\in V(H^{i})$ with degree
less than $k$, consider an edge $f$ that consists of all $k$ copies
of $v$ in $H^{*}$. Add all these new edges to $H^{*}$ and name
obtained hypergraph $H^{i+1}$. Note that $k-\delta(H^{i+1})=t-1$
and therefore $H'$ will be constructed in finite steps. Therefore,
any $k$-uniform hypergraph $H$ ($k\geq4$) that has the maximum
degree at most $k$, is 2-colorable. Note that for any graph $G$,
$\chi_d(G)\leq \Delta(G)+1$ unless $G=C_5$, in which case
$\chi_d(C_5)=5$ (see \cite{MR1991048}). Regarding the above
discussion and the proof of Lemma \ref{delta}, for any
$k$-regular graph $G$, $\chi_d(G)\leq 2\chi(G)$. This result and
Theorem \ref{main} imply the next corollary.
\begin{cor}
There exists a constant $c$ such that for any $k$-regular graph
$G$, $\chi_d(G)\leq \min\{\ 2\chi(G),\ \chi(G)+c\ln k\}$.
\end{cor}
In the rest of this paper, we are focused on relationship between the
 total dominating set (res. double total dominating set) and the dynamic chrommatic number
of graphs.
\begin{lem}\label{domin}
Let $G$ be a graph.
\begin{enumerate}
\item If  $e\Delta^2\leq2^{\delta(G)-1}$ and there is a total dominating set
    $T\subsetneq V(G)$ then
    $\chi_d(G)\leq \chi(G[V\setminus T])+2\chi(G[T])$.
\item If $G$ has a  double total dominating set $T\subseteq V(G)$ then $\chi_d(G)\leq\chi(G[V\setminus T])+\chi(G[T])$.
\end{enumerate}
\end{lem}
\begin{proof}{
 For convenience let $\chi(G[T])=s$.

(1) Assume that $H$ is a hypergraph with the
vertex set $T$ and the hyperedge set defined as follows.
$$E(H)\isdef\{N(u)|\ N(u)\subseteq T,\ u\in V(G)\}.$$  By proof of
Lemma \ref{delta}, $H$ is $2$-colorable. Assume that $c$ is an
$s$-coloring of $G[T]$ and $f$ is a $2$-coloring of $H$.
Obviously, $h'\isdef(c,f)$ is a $2s$ coloring of $G[T]$.
Consider the coloring $h$ of $G$ in which the restriction of $h$
to $T$ is the same as $h'$ and the vertices in $G[V\setminus T]$
are colored with the colors that are not used in $h'$. One can
easily check that the coloring $h$ is a dynamic coloring of $G$.

(2) Let $h'$ be a coloring of $G[T]$ and $h$ be a coloring that is the same on
$T$ as $h'$  and the vertices in $V\setminus T$ are colored by
$h$ with colors that are distinct from the colors used in
$T$. It is easy to see that $h$ is a dynamic coloring of $G$ and
therefore $\chi_d(G)\leq \chi(G[V\setminus T])+\chi(G[T])$.
}\end{proof}
As a consequence of the previous lemma, when a graph $G$ has a double total dominating set $T$,
 $G$ does not have a dynamic chromatic
number far from its chromatic number. This result is restated in
the next corollary.
\begin{cor}\label{domincor}
Let $G$ be a graph and assume that there is a double total
dominating set $T\subsetneq V(G)$ . Then $\chi_d(G)\leq
2\chi(G)$.
\end{cor}
For a Kneser graph $KG(m,n)$, if $m\geq 3n$ then every vertex of
$KG(m,n)$ is in some triangles and therefore
$\chi_d(KG(m,n))=\chi(KG(m,n))$. But when $m=2n+t<3n$ the graph
$KG(m,n)$ is triangle free and so the dynamic chromatic number of
$KG(m,n)$ is still interesting.
Let $T=\{A\in V(KG(m,n))| A\subseteq \{t+1,t+2,\ldots,m\}\}$.
Obviously, the induced subgraph  $G[T]$ is a bipartite
graph and $T$ is a total dominating
set. Note that $c(B)\isdef \min B$ is a $t$-coloring of
the induced subgraph on $V(KG(m,n))\setminus T$.
By Lemma \ref{domin}, we have $\chi_d(KG(m,n))\leq
t+4=\chi(KG(m,n))+2$. Although, the exact value of
$\chi_d(KG(m,n))$ is not determined,
but it is shown that the Conjecture
\ref{conj} is true for Kneser graphs.

\begin{cor}
Let $G$ be a $k$-critical graph. If $e\Delta(G)^2\leq
2^{\delta(G)-1}$ then $\chi_d(G)\leq 2k-2$
\end{cor}
\begin{proof}{
Let $T$ be a random subset of $V(G)$ such that each vertex $v$
lies  in $T$  with probability ${1\over 2}$, randomly and
independently. Let ${\cal A}_x$ be the event that {\it all
neighbors of $x$ are in $T$ or none of them are in $T$}. For each
vertex $v\in V(G)$, ${\cal A}_v$ is mutually independent of all
but at most $\Delta(G)(\Delta(G)-1)$ events. Lovasz Local Lemma
guarantees that with positive probability $\bigcap \bar{{\cal
A}_x}$ happens. Equivalently, there is a $T\subsetneq V(G)$ such that
for every vertex $u\in V(G)$, ${\rm deg}_T(u)>0$ and ${\rm
deg}_{V(G)\setminus T}(u)>0$. Since $G$ is a $k$-critical graph,
$\chi(G[T])< \chi(G)$ and also $\chi(G[V(G)\setminus T])<
\chi(G)$. By second part of Lemma \ref{domin}, the proof is
completed. }
\end{proof}

In the proof of the previous corollary it is shown that if for a
graph $G$, $e\Delta(G)^2\leq 2^{\delta(G)-1}$ then $G$ has a total
dominating set. Therefore, by Lemma \ref{domincor}, $\chi_d(G)\leq
2\chi(G)$. This provides another proof of Lemma \ref{delta}.\\

\noindent{\bf Acknowledgment}\\

The author would like to thank Hossein Hajiabolhassan and Saeed Shaebani for their
invaluable comments.


\end{document}